\documentclass[10pt,letterpaper,conference]{IEEEtran}
\makeatletter
\def\bstctlcite{\@ifnextchar[{\@bstctlcite}{\@bstctlcite[@auxout]}}
\def\@bstctlcite[#1]#2{\@bsphack
  \@for\@citeb:=#2\do{%
    \edef\@citeb{\expandafter\@firstofone\@citeb}%
    \if@filesw\immediate\write\csname #1\endcsname{\string\citation{\@citeb}}\fi}%
  \@esphack}
\makeatother

\usepackage{cite}

\ifCLASSINFOpdf
\usepackage[pdftex]{graphicx}
\graphicspath{{./pdf/}{./jpeg/}}
\DeclareGraphicsExtensions{.pdf,.jpeg,.png}
\else
\usepackage[dvips]{graphicx}
\graphicspath{{./eps/}}
\DeclareGraphicsExtensions{.eps}
\fi
\usepackage{amsmath}
\usepackage{amsfonts}
\usepackage{multirow}
\usepackage{hhline}
\usepackage{xcolor,colortbl}
\definecolor{Gray}{gray}{0.85}
\usepackage{array}

\usepackage{stfloats}
\usepackage{url}

\IEEEoverridecommandlockouts

\begin{document}
\bstctlcite{IEEEexample:BSTcontrol}
%
\title{Optimal Selection of Small-Scale Hybrid PV-Battery Systems to Maximize Economic Benefit Based on Temporal Load Data}

\author{\IEEEauthorblockN{Jeremy Every, Li Li}
\thanks{This research is supported by an Australian Government Research Training Program Scholarship}
\IEEEauthorblockA{Faculty of Engineering and Information Technology\\
University of Technology Sydney\\
Ultimo 2007, Australia \\
Email: jeremy.every@student.uts.edu.au}
\and
\IEEEauthorblockN{David G. Dorrell}
\IEEEauthorblockA{College of Agriculture, Engineering and Science\\
University of KwaZulu-Natal\\Durban 4041, South Africa\\
Email: dorrelld@ukzn.ac.za}}


%


\maketitle

\begin{abstract}
Continued advances in PV and battery energy storage technologies have made hybrid PV-battery systems an attractive prospect for residential energy consumers. However the process to select an appropriate system is complicated by the relatively high cost of batteries, a multitude of available retail electricity plans and the removal of PV installation incentive schemes. In this paper, an optimization strategy based on an individual customer's temporal load profile is established to maximize electricity cost savings through optimal selection of PV-battery system size, orientation and retail electricity plan. Quantum-behaved particle swarm optimization is applied as the underlying algorithm given its well-suited application to problems involving hybrid energy system specification. The optimization strategy is tested using real-world residential consumption data, current system pricing and available retail electricity plans to establish the efficacy of a hybrid PV-battery solution.
\end{abstract}


%
\IEEEpeerreviewmaketitle

\section{Introduction}
Energy storage systems for residential applications, particularly lithium-ion based battery systems, have undergone rapid development in the past five years. A range of stationary energy storage devices from manufacturers including Tesla Motors, Enphase, Mercedez Benz, Samsung and LG have been recently introduced to the market. In Australia, average annual electricity prices have risen 4.5\% over the last 10 years \cite{ABS16}, 2.5\% higher than the inflation rate target \cite{RBA16}. Consequently the application of battery systems to complement existing rooftop PV systems or the installation of new hybrid PV-battery solutions are of particularly interest to energy consumers aiming to reduce their net electricity bills.

The performance of PV-battery systems in reducing energy costs have been investigated in \cite{Ren16,Parra16,Mulder13,Beck16} where prescribed PV and battery sizes were tested under various tariffs \cite{Ren16,Parra16}, incentive schemes \cite{Mulder13} and temporal resolution of energy consumption data \cite{Beck16}. However existing research has been primarily structured around typical PV-battery systems and consumption profiles, rather than investigating strategies to enable individual customers to optimize a system based on their own circumstances.

In this paper an optimization methodology is developed with the aim to assist potential PV-battery investors in determining the economic efficacy of a hybrid PV-battery system. Models for solar insolation, PV energy output, battery operation, installation cost and maintenance cost are presented to form key components of the underlying optimization problem. In order to maximize the net benefit of a hybrid system, the PV power rating, PV orientation, retail electricity plan and currently available retail battery products and battery operating modes are optimally selected. The objective function is formulated as a net present value (NPV) evaluation of the electricity cost savings that can be achieved through the introduction of a hybrid PV-battery system compared to a known lowest cost retail electricity plan. 
 
A form of particle swarm optimization (PSO), known as quantum-behaved PSO (QPSO), is utilized due to its fast convergence speed, natural handling of optimization parameter constraints and simplicity of implementation. The algorithm is tested in an Australian context for three residences in the state of New South Wales using one year of hourly energy data and five years of solar insolation data from the Australian Bureau of Meteorology. Time-of-use (TOU) retail electricity plans from three large Australian retailers as well as current PV and retail battery pricing are applied in order to establish the feasibility of installing a PV-battery system under prevailing market conditions. 

\section{Model Definition}

\subsection{Solar Insolation Model}
\label{ss:solar}
The total hourly insolation $I_{T}$ on a tilted plane includes contributions from beam insolation, diffuse insolation and a ground reflected component. A multitude of transposition models have been developed enabling the estimation of insolation on a tilted plane based on horizontal insolation data. Among the available models, the Hay-Davies-Klucher-Reindl (HDKR) model, represented by (\ref{eq:tiltinsol}), has been identified as one of the more accurate \cite{Noorian08,Duffie13}.
\begin{IEEEeqnarray}{Rl}
I_{T}=&\left(I_{b}+A_{i}I_{d}\right)R_{b}\nonumber\\
&+\:I_d(1-A_{i})\left(\frac{1+\cos\beta}{2}\right)\left[1+f\sin^{3}\left(\frac{\beta}{2}\right)\right]\nonumber\\
&+\:I\rho_{g}\left(\frac{1-\cos\beta}{2}\right)\label{eq:tiltinsol}%
\end{IEEEeqnarray}

In (\ref{eq:tiltinsol}), $I_{b}$ and $I_{d}$ are the hourly beam and diffuse insolation on a horizontal plane respectively, $A_i=I_b/I_o$, $f=\sqrt{I_b/I}$, $I$ is the global horizontal radiation, $I_o$ is the hourly extraterrestrial insolation incident on a horizontal plane projected from the Earth's surface, $\rho_g$ is the ground reflectance and $R_b$ is the ratio of tilted to horizontal beam radiation. Importantly, $R_b$ is a function of panel tilt $\beta$ and panel azimuth $\gamma$, equations for which are well established in literature and can be found in \cite{Duffie13}. Consequently tilt and azimuth must be optimized in order to maximize the incident insolation $I_T$ on a PV array. 

Finally, it should be noted that $\gamma = 0^\circ$ implies a north facing surface in this paper.

\subsection{Photovoltaic Model}
\label{ss:photo}

According to Duffie and Beckman \cite{Duffie13}, the output energy of a PV system is defined as:
\begin{equation}\label{eq:pvenergy}
E_{pv}=A_{c}ZI_{T}\eta_{mpp}\eta_e\zeta_{pv}
\end{equation}
where $A_c$ is the PV panel area, $Z$ in the number of panels, $\eta_{mpp}$ is the PV panel operating efficiency, $\eta_e$ is the efficiency of the associated balance of plant and $\zeta_{pv}$ is the annual degradation factor of the PV panels. Although not explicitly defined in this paper, $\zeta_{pv}$ is assumed to be linear, in-line with the manufacturer's warranty. 

The operating efficiency $\eta_{mpp}$ is defined as:
\begin{equation}\label{eq:opeff}
\eta_{mpp}=\eta_{mpp,STC}+\mu_{mpp}(T_{c}-T_{a})
\end{equation}
where $T_a$ and $T_c$ are the ambient and cell temperatures respectively, $\eta_{mpp,STC}$ is the panel efficiency at standard test conditions and $\mu_{mpp}$ (\%/W) is the power coefficient. Both $\eta_{mpp,STC}$ and $\mu_{mpp}$ are defined in manufacturer data sheets. 

The cell temperature $T_c$ is defined as:
\begin{equation}\label{eq:celltemp}
T_{c}=T_{a}+(T_{NOCT}-20)\cdot\frac{G_{T}}{800}\cdot(1-\eta_{mpp,STC})
\end{equation}
where $G_{T}$ is the incident irradiance (assumed to be uniform over an hour period and therefore equal to $I_{T}$) and $T_{NOCT}$ is the cell temperature under nominal operating cell temperature (NOCT) conditions as detailed in manufacturer data sheets.

In an Australian context, PV system costs are subsidized through small-scale technology certificates (STCs). The quantity of STCs generated after installation is based on the power rating $P_{pv,rat}$ and a location multiplier $M_{loc}$. In this paper $M_{loc}$ is assumed to be 20.73 (for the east coast of New South Wales) while STCs are assumed to be worth $C_{STC}$ = \$34 \cite{Jacobs16}. Consequently in Australia, the net system cost $S_{pv}$ is defined as:
\begin{equation}\label{eq:pvsystemcost}
S_{pv}=U_{pv}P_{pv,rat}-M_{loc}P_{pv,rat}C_{STC}
\end{equation}
where $U_{pv}$ is the average price per watt peak. Based on data provided in \cite{Jacobs16}, March 2016 $U_{pv}$ prices were (in AUD) \$3.20, \$3.00, \$2.55, \$2.35 and \$2.20 for 1 kW, 1.5 kW, 3 kW, 5~kW and 10 kW rated systems respectively. When solving the optimization problem presented in Section~\ref{s:optim}, the price corresponding to the closest system size is used.

The PV modules considered in this research were modelled based on 280 W Trina Solar TSM-PC05A polycrystalline modules.

\subsection{Battery Model}
\label{ss:battery}

The battery models defined in this section are structured based on manufacturer warranties and guarantees to establish the economic benefit the owner can expect over the lifetime of the system.

The maximum capacity of a battery decreases over its lifetime due to a number of factors, chief of which is the number of charge/discharge cycles undergone. The degradation rate $\zeta_{batt}$ (kWh/cycle) is defined as:
\begin{IEEEeqnarray}{Rl}
\zeta_{batt}=(C_{\max 0}-C_{EOL})/Y_{EOL}\label{eq:battdegrad}
\end{IEEEeqnarray}
where $C_{\max 0}$, $C_{EOL}$ and $Y_{EOL}$ are the initial maximum capacity, end-of-life maximum capacity and cycle life respectively as defined in manufacturer data sheets. 

The maximum capacity $C_{\max,qdh}$ available at the start of each hour $h$, in day $d$ and billing period $q$, is assumed to be a linear function of the number of operational cycles $Y_{qdh}$ in the previous hour such that:
\begin{IEEEeqnarray}{Rl}
C_{\max,qdh}=C_{\max,qd(h-1)}-Y_{qd(h-1)}\zeta_{batt}\label{eq:maxcap}
\end{IEEEeqnarray}

It should be noted that $Y_{qdh}$ generally represents only a partial cycle for hourly intervals and therefore represents only a fraction of the energy throughput of a full discharge/charge cycle. $Y_{qdh}$ is defined as follows:
\begin{IEEEeqnarray}{Rl}
Y_{qdh}=\frac{E_{bpv,qdh}+E_{bg,qdh}+E_{bd,qdh}}{2DC_{\max,qdh}}\label{eq:cycles}
\end{IEEEeqnarray}
where $E_{bd,qdh}$, $E_{bpv,qdh}$ and $E_{bg,qdh}$ are the discharge, PV charge and grid charge energy flows respectively and $D$ is the maximum depth of discharge.
 
The available capacity at the start of each hour is a function of the capacity at the start of the previous hour and the total charge/discharge energy that has flowed to/from the battery cells in the previous hour. The available capacity $C_{qdh}$ is therefore defined as:
\begin{IEEEeqnarray}{Rl}
C_{qdh}=&C_{qd(h-1)}-E_{bd,qd(h-1)}+E_{bpv,qd(h-1)}\nonumber\\
&{}+E_{bg,qd(h-1)}\label{eq:capacity}
\end{IEEEeqnarray}

The charge and discharge energy flow terms are defined in (\ref{eq:pvcharge})\textendash(\ref{eq:discharge}).
\begin{IEEEeqnarray}{Rl}
E_{bpv,qdh}=\max\Bigl\lbrace \min\bigl[&C_{\max,qdh}-C_{qdh},\nonumber\\
&{}(E_{pv,qdh}-E_{load,qdh})(1-F),\nonumber\\
&{}R_{\max}(1-F)\bigr],0\Bigr\rbrace\label{eq:pvcharge}
\end{IEEEeqnarray}
\begin{IEEEeqnarray}{Rl}
E_{bg,qdh}=\max\Bigl\lbrace \min\bigl[&C_{\max,qdh}-C_{qdh},\nonumber\\
&{}R_{\max}(1-F)\bigr](M_3+M_4)I_{op,qdh}\nonumber\\
&{}-E_{bpv,qdh},0\Bigr\rbrace\label{eq:gridcharge}
\end{IEEEeqnarray}
\begin{IEEEeqnarray}{Rl}
E_{bd,qdh}=\max\Bigl\lbrace \min\bigl[&C_{qdh}-C_{\max,qdh}(1-D),\nonumber\\
&{}(E_{load,qdh}-E_{pv,qdh})/(1-F),R_{\max}\bigr]\nonumber\\
&{}\times\bigl[(M_2+M_4)I_{sh,qdh}\nonumber\\
&{}\quad\;\;+I_{pk,qdh}\bigr],0\Bigr\rbrace\label{eq:discharge}
\end{IEEEeqnarray}

In (\ref{eq:pvcharge})\textendash(\ref{eq:discharge}), $E_{pv,qdh}$, $E_{load,qdh}$ and $R_{\max}$ are the PV generated energy as defined in (\ref{eq:pvenergy}), load energy demand and rated continuous charge/discharge rate respectively. $I_{op,qdh}$, $I_{sh,qdh}$, $I_{pk,qdh}$ and terms of the form $M_x$ are battery operation control variables defined later.

It should be noted that the charge energy flow terms $E_{bpv,qdh}$ and $E_{bg,qdh}$ are considered to be the net additional charge to a battery after losses while the discharge energy term $E_{bd,qdh}$ is the total energy discharged from the battery (i.e. usable energy plus losses). The losses have been accounted for through the inclusion of a loss factor $F=(1-\eta_{batt})/2$ where $\eta_{batt}$ is the battery round-trip efficiency.

The total energy loss due to battery charging and discharging is:
\begin{equation}\label{eq:bloss}
E_{bloss,qdh}=E_{bpvloss,qdh}+E_{bgloss,qdh}+E_{bdloss,qdh}
\end{equation}
where $E_{bpvloss,qdh}$, $E_{bgloss,qdh}$ and $E_{bdloss,qdh}$ are the losses during PV charging, grid charging and discharging respectively defined in (\ref{eq:pvloss})\textendash(\ref{eq:dischargeloss}).
\begin{IEEEeqnarray}{Rl}
E_{bpvloss,qdh}=\max\Bigl\lbrace\min\bigl[ & (C_{\max,qdh}-C_{qdh})/(1-F),\nonumber\\
&{}E_{pv,qdh}-E_{load,qdh},\nonumber\\
&{}R_{\max}\bigr],0\Bigr\rbrace F\label{eq:pvloss}
\end{IEEEeqnarray}
\begin{IEEEeqnarray}{Rl}
E_{bgloss,qdh}=\max\Bigl\lbrace\min\bigl[&(C_{\max,qdh}-C_{qdh})/(1-F),R_{\max}\bigr]\nonumber\\
&{}\times (M_3+M_4)I_{op,qdh}\nonumber\\
&{}-E_{bpv,qdh},0\Bigr\rbrace F\label{eq:gridloss}
\end{IEEEeqnarray}
\begin{IEEEeqnarray}{Rl}
E_{bdloss,qdh}=\max\Bigl\lbrace \min\bigl[&E_{load,qdh}-E_{pv,qdh},R_{\max},\nonumber\\
&{}C_{qdh}-C_{\max,qdh}(1-D)\bigr]\nonumber\\
&{}\times\bigl[(M_2+M_4)I_{sh,qdh}\nonumber\\
&{}\quad\;\;+I_{pk,qdh}\bigr],0\Bigr\rbrace\ F\label{eq:dischargeloss}
\end{IEEEeqnarray}

In addition to increasing the self-consumption ratio of PV generated energy (by effectively translating PV generated energy to non-generation periods), an energy storage system can also be used to perform energy arbitrage by charging during low cost off-peak hours and discharging during peak periods. In this paper, a review of various battery operating modes is undertaken to determine the most economically efficient mode for each residence assessed. The operating modes considered are defined as follows:
\begin{itemize}[\IEEEsetlabelwidth{Mode 4:}]
\item[Mode 1:] PV generation shifting. Discharge in peak only.
\item[Mode 2:] PV generation shifting. Discharge during shoulder and peak periods.
\item[Mode 3:] Energy arbitrage and PV generation shifting. Discharge in peak only.
\item[Mode 4:] Energy arbitrage and PV generation shifting. Discharge during shoulder and peak periods.
\end{itemize}

As previously indicated, (\ref{eq:gridcharge}), (\ref{eq:discharge}), (\ref{eq:gridloss}) and (\ref{eq:dischargeloss}) are controlled by the operation mode variables $M_1$, $M_2$, $M_3$ and $M_4$ where:
\begin{equation}\label{eq:mode}
M_x=\begin{cases}
& 1 \quad \text{if in Mode}\;x\\
& 0 \quad \text{otherwise}
\end{cases}
\end{equation}

The variables $I_{op,qdh}$, $I_{sh,qdh}$ and $I_{pk,qdh}$ control battery charge and discharge based on the tariff period within which a particular hour lies and take the form:
\begin{equation}\label{eq:offpeak}
I_{op,qdh}=\begin{cases}
& 1 \quad \text{if}\;\;h\in\lbrace\text{offpeak hours}\rbrace\\
& 0 \quad \text{otherwise}
\end{cases}
\end{equation}
with similar equations for $I_{sh,qdh}$ and $I_{pk,qdh}$ for shoulder and peak hours respectively.

The final component of the battery model is the battery cost, defined simply as:
\begin{equation}\label{eq:battsystemcost}
S_{b}=U_{b}X
\end{equation}
where $U_b$ is the price per battery and $X$ is the number of battery units installed.

Two battery systems were considered in this research \textendash{ }the Tesla Motors 13.5 kWh 5 kW Powerwall 2 and the more modular 1.2 kWh 260 W Enphase AC Battery. Values for the battery model parameters defined in this section including $D$, $Y_{life}$, $Y_{EOL}$, $C_{EOL}$ and $\eta_{batt}$ were based on the manufacturer datasheets \cite{TeslaPWDC16} and \cite{Enphase16} for the Tesla and Enphase systems respectively. Based on pricing provided by Tesla, the fully installed cost of the Powerwall 2 is AU\$10,000 \cite{TeslaPW16}, while the cost of the Enphase battery is approximately AU\$2,000 \cite{SQ16}.

\subsection{Maintenance Model}
\label{ss:maitenance}

During the lifetime of a PV-battery system, periodic maintenance as well as battery and inverter replacements are required. In this paper, the lifespans of PV modules and inverters/batteries are assumed to be 20 years and 10 years respectively. Consequently with $t$ billing periods per year, inverters/batteries will require replacement after $10t$ billing periods. Furthermore periodic maintenance is assumed to occur every $5t$ billing periods. The system maintenance costs are therefore defined as:
\begin{equation}\label{eq:maintcost}
W_q=\begin{cases}
200\;\;\text{if}\;\frac{q-1}{5t}\in\mathbb{Z}^+\text{, }\frac{q-1}{10t}\notin\mathbb{Z}^+\\
400+\kappa_{inv}U_{inv}S_{pv}+\kappa_bS_{b}\;\;\text{if}\;\frac{q-1}{10t}\in\mathbb{Z}^+\\
0\;\;\text{otherwise}
\end{cases}
\end{equation}
where $U_{inv}$ is the inverter replacement cost (\$/W$_\text{ac}$) and $\kappa_{inv}$ and $\kappa_b$ are cost reduction factors for the inverter and batteries. Given a current average per unit inverter cost of US\$0.29/W$_\text{ac}$ \cite{NREL15}, the per unit PV inverter costs are assumed to be $U_{inv}=\text{AU}\$0.41$/W$_\text{ac}$ (assuming AUD/USD exchange rate of 0.7). The cost of inverters and batteries are forecasted to reduce significantly over the next 10 years with reductions of 31\% and 53\% respectively between 2015-2025 \cite{CSIRO15}. Consequently the cost reduction factors in (\ref{eq:maintcost}) are assumed to be $\kappa_{inv}=0.69$ and $\kappa_b=0.47$.

\section{Optimization Problem}
\label{s:optim}

The objective of this research is to maximize the electricity cost savings achieved through optimal selection of a hybrid PV-battery system based on high resolution smart meter load data and prevailing economic and PV-battery market conditions. Cost savings are quantified through an NPV analysis performed on the difference in electricity costs between a known lowest cost retail electricity plan and a hybrid PV-battery system combined with other currently available retail electricity plans. 

As previously indicated in Section~\ref{ss:battery} and \ref{ss:maitenance}, hourly evaluations of the energy flows are conducted for each hour $h$ in day $d$ and billing period $q$. Maximizing the net benefit over all the billing periods $Q$ in the lifetime of the system is the objective of the optimization problem as defined below.

\subsection{Problem Definition}
\label{ss:probdef}

\vspace*{+0.5\baselineskip}
\noindent\emph{Given:}
\begin{enumerate}
\item Latitude and longitude of the location
\item Hourly load and insolation profile
\item A real annual discount rate of 3.92\% (6\% nominal rate and 2\% inflation)
\item Real annual electricity price growth of 2\%
\item PV system balance of plant efficiency (90\%)
\item System lifespan (20 years)
\end{enumerate}

\vspace*{+0.5\baselineskip}
\noindent\emph{Find:} Tilt angle $\beta$, azimuth angle $\gamma$, number of PV panels $Z$ and number of batteries $X$

\vspace*{+0.5\baselineskip}
\noindent\emph{Objective:}
\begin{IEEEeqnarray}{Rl}
\max_{\beta,\gamma,Z,X}NPV=&\sum_{q=1}^{Q}\frac{\left(C_{base,q}-C_{pvbatt,q}\right)\left( 1+r_{e}\right)^q}{\left(1+r_{d}\right)^{q}}\nonumber\\
&-\:\sum_{q=1}^{Q}\frac{W_{q}}{\left(1+r_{d}\right)^{q}}-\bigl( S_{pv}+S_{b}\bigr)\label{eq:objfunc}%
\end{IEEEeqnarray}

\noindent\emph{Subject to:}
\begin{IEEEeqnarray}{rCll}
0 \leq & \beta & \leq 180 &\quad\text{for}\;\beta\in\mathbb{R}\IEEEyesnumber\IEEEyessubnumber\label{eq:tiltcons}\\
-180 < & \gamma & \leq 180 &\quad\text{for}\;\gamma\in\mathbb{R}\IEEEyessubnumber\label{eq:azicons}\\
0 \leq & Z & \leq Z_{\max} &\quad\text{for}\;Z\in\mathbb{Z}^{+}\IEEEyessubnumber\label{eq:panelcons}\\
0 \leq & X & \leq X_{\max} &\quad\text{for}\;X\in\mathbb{Z}^{+}\IEEEyessubnumber\label{eq:battcons}%
\end{IEEEeqnarray}

In (\ref{eq:objfunc}), $C_{pvbatt,q}$ and $C_{base,q}$ are the cost of electricity with and without a PV-battery system within the billing period $q$. As quarterly billing is assumed in this paper, the real discount rate 3.92\% and the real annual electricity price growth of 2\% are adjusted to the quarterly effective rates $r_{d}=0.97\%$ and $r_{e}=0.50\%$ respectively. 

The terms $C_{base,q}$ and $C_{pvbatt,q}$ are defined as:
\begin{equation}\label{eq:basecost}
C_{base,q}=\sum_{d=1}^{D}\left( \sum_{h=1}^{24}T_{grid0,qdh}E_{load,qdh}+T_{sc0,qd}\right)
\end{equation}
\begin{IEEEeqnarray}{Rl}
C_{pvbatt,q}=\sum_{d=1}^{D}\Bigg\lbrace\sum_{h=1}^{24}&\Big[T_{grid,qdh}\max\left(0,E_{bal,qdh}\right)\nonumber\\
&{}-T_{feed,qdh}\max\left(0,-E_{bal,qdh}\right)\Big]\nonumber\\
&{}+T_{sc,qd}\Bigg\rbrace\label{eq:pvbattcost}
\end{IEEEeqnarray}
where $T_{grid0,qdh}$ and $T_{grid,qdh}$ are the grid imported electricity tariff of the base plan and tested plan respectively for the $h^{th}$ hour of day $d$ with $D$ days in the billing period, $T_{sc0,qd}$ and $T_{sc,qd}$ are the daily electricity supply charges for the base plan and tested plan respectively, $T_{feed,qdh}$ is the PV feed-in tariff and $E_{bal,qdh}$ is the net energy flow balance of the terms defined in Section~\ref{ss:battery}, expressed as:
\begin{IEEEeqnarray}{Rl}
E_{bal,qdh}=&E_{load,qdh}-E_{pv,qdh}-E_{bd,qdh}\nonumber\\
&{}+E_{bpv,qdh}+E_{bg,qdh}+E_{bloss,qdh}\label{eq:energybal}
\end{IEEEeqnarray}

As the optimization parameters $Z$ and $X$ are limited to integer values (with respective maximums $Z_{\max}$ and $X_{\max}$ determined by the customer's available space restrictions), while $\beta$ and $\gamma$ may take any real value within the domain of the constraints, the problem is classified as a mixed-integer non-linear programming (MINLP) problem.  

\subsection{Optimization Algorithm}
\label{ss:optalg}

To solve MINLP problems, metaheuristic programming methods such a PSO have seen increased application in the last decade, with PSO already applied to a number of PV optimization problems \cite{Yadav13}. PSO simulates the social interaction within bird flocks to achieve a global objective without a governing central controller \cite{Sun12}. A modified version of PSO described by Sun et al. \cite{Sun12} known as Quantum-behaved PSO (QPSO), uses the principles of quantum mechanics, in particular quantum delta potential wells, to sample around the best positions to eventually find the global best position. The defining algorithm for QPSO is relatively simple and requires fewer parameter adjustments between problems than PSO \cite{Sun12}. Furthermore, QPSO handles parameter constraints naturally with no specific modifications required to the algorithm, as opposed to PSO.

In order to handle the discrete parameters $Z$ and $X$, the hypercube nearest-vertex approach adopted in \cite{Chowdhury13} was utilized, which effectively rounds the position of the candidate particle to the nearest integer prior to evaluating the objective function. For installation simplicity, the tilt and azimuth angles were also considered as discrete in this analysis. 

The QPSO optimization algorithm was developed and simulated in Matlab version R2015b.

\begin{table*}[!hb]
\renewcommand{\arraystretch}{1.3}
\caption{Characteristics and Economic Performance of Optimized PV-Battery System For Different Retail Electricity Plans}
\label{tab:results}
\centering
\begin{tabular}{|c|c|c|c|c|c|c|c|c|}
\hline
Customer & Retailer & Battery Size (kWh) & PV Size (kW$_\text{p}$) & Tilt & Azimuth & NPV  & MIRR & Payback (Years) \\\hline\hline
\multirow{3}{*}{1} & A & 0 & 8.4 & $29^\circ$ & $30^\circ$ & \$10,534 & 7.39\% & 9.3\\\hhline{~*{8}{-}}
  & \cellcolor[gray]{0.8}B & \cellcolor[gray]{0.8}0 & \cellcolor[gray]{0.8}8.4 & \cellcolor[gray]{0.8}$29^\circ$ & \cellcolor[gray]{0.8}$30^\circ$ & \cellcolor[gray]{0.8}\$11,153 & \cellcolor[gray]{0.8}7.52\% & \cellcolor[gray]{0.8}9.0\\\hhline{~*{8}{-}}
 & C & 0 & 8.4 & $29^\circ$ & $30^\circ$ & \$10,745 & 7.44\% & 9.0\\\hline\hline

\multirow{3}{*}{2} & A & 0 & 4.2 & $31^\circ$ & $26^\circ$ & \$672 & 4.88\% & 18.4\\\hhline{~*{8}{-}}
  & \cellcolor[gray]{0.8}B & \cellcolor[gray]{0.8}0 & \cellcolor[gray]{0.8}4.2 & \cellcolor[gray]{0.8}$31^\circ$ & \cellcolor[gray]{0.8}$26^\circ$ & \cellcolor[gray]{0.8}\$917 & \cellcolor[gray]{0.8}5.04\% & \cellcolor[gray]{0.8}17.8\\\hhline{~*{8}{-}}
 & C & 0 & 4.2 & $31^\circ$ & $26^\circ$ & \$732 & 4.92\% & 18.2\\\hline\hline

\multirow{3}{*}{3} & A & 0 & 7.56 & $29^\circ$ & $25^\circ$ & \$4,601 & 6.06\% & 13.8\\\hhline{~*{8}{-}}
  & \cellcolor[gray]{0.8}B & \cellcolor[gray]{0.8}0 & \cellcolor[gray]{0.8}7.56 & \cellcolor[gray]{0.8}$29^\circ$ & \cellcolor[gray]{0.8}$25^\circ$ & \cellcolor[gray]{0.8}\$5,044 & \cellcolor[gray]{0.8}6.19\% & \cellcolor[gray]{0.8}13.6\\\hhline{~*{8}{-}}
 & C & 0 & 7.56 & $30^\circ$ & $26^\circ$ & \$4,870 & 6.14\% & 13.6\\\hline
\end{tabular}
\end{table*}

\section{Input Data}
\label{s:input}

Electricity consumption data measured over a one-year period for three arbitrarily selected customers from the Australian Government initiated `Smart Grid, Smart City' project \cite{Arup14} was used in the analysis. Daily insolation and ambient temperature data over a five year period for each location were derived from the Australian Bureau of Meteorology Climate Data Online database \cite{BOM16}. The daily data was converted to hourly data using the methodology established in \cite{Every16}. 

The electricity tariff structures tested in the optimization problem were based on real  2016 TOU rates from three large Australian retailers.

\section{Results and Discussion}
\label{s:results}

Table~\ref{tab:results} shows a summary of the optimized PV-battery systems for each customer under retail electricity plans from three large retailers. For each customer, Retailer B was found to provide the greatest benefit. For Customer 1, the maximum system size of 8.4 kW ($Z_{\max}=30$, 280 W PV panels) was reached, while Customer 3 was also found to potentially benefit from a relatively large PV system of 7.56 kW. In contrast, the optimal PV system for Customer 2 was found to be a far smaller system at 4.2 kW. The optimal tilt and azimuth angles also varied for each customer but were found to be in 20-30 degree range. However importantly for the customers assessed, no instance was found whereby an energy storage system would yield an economic benefit higher than a PV-only system based on current battery pricing. 
 
\begin{figure}[!t]
\centering
\includegraphics[width=\linewidth]{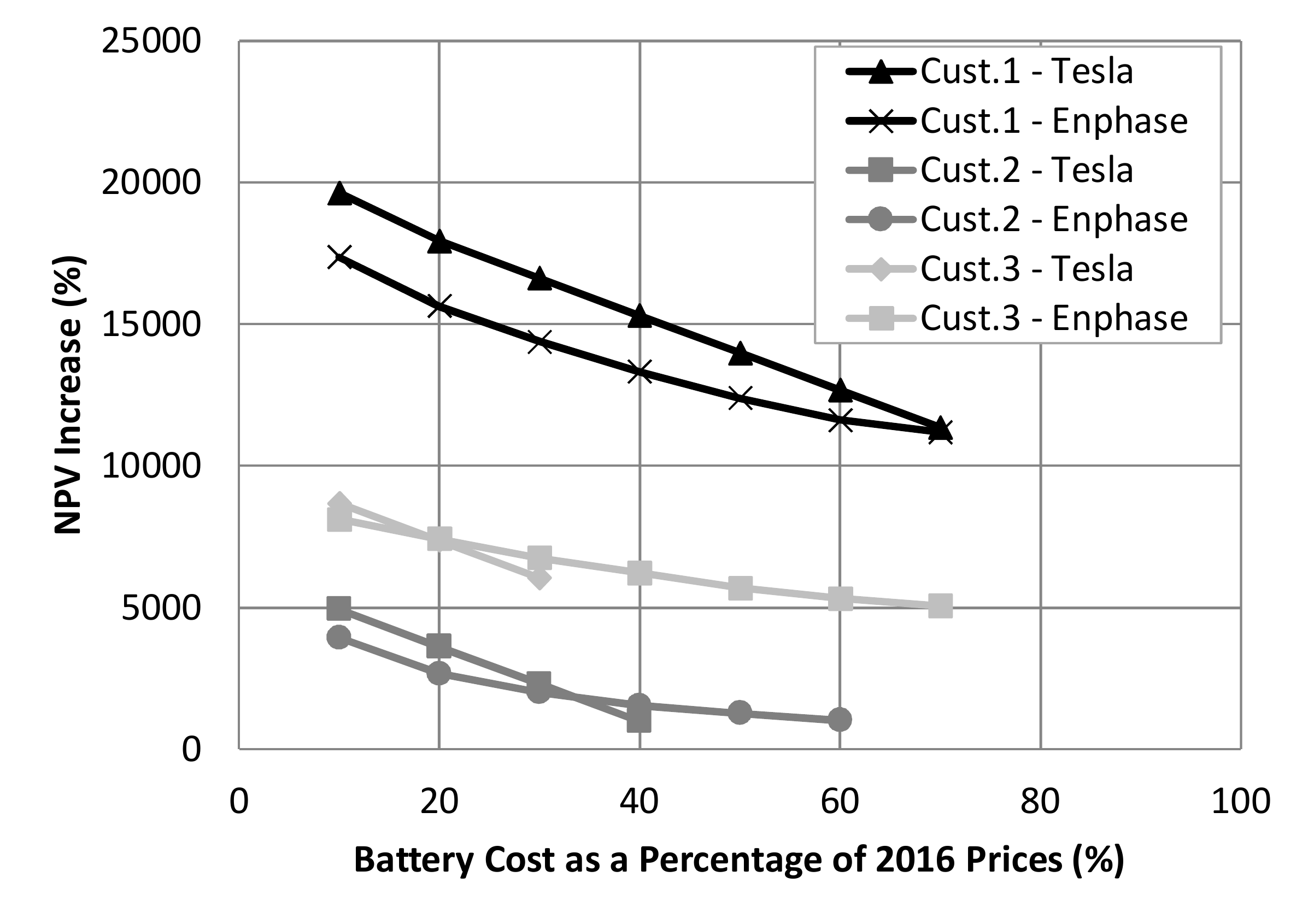}
\caption{NPV sensitivity to installed battery cost (Retailer B, Mode 2 operation)}
\label{fig:npvsensi}
\end{figure}

To determine the price point at which a hybrid PV-battery system becomes an economically beneficial option for each customer, an NPV sensitivity analysis was undertaken on battery pricing. Referring to Fig.~\ref{fig:npvsensi}, systems consisting of either the Tesla Powerwall 2 or an Enphase AC battery become viable for Customer~1 when unit costs are reduced to 70\% of 2016 pricing. Customers 2 and 3 would first see a benefit from the small, more modular Enphase system at the 60-70\% price point but would not see a benefit from a larger Tesla system until pricing reached 30-40\% of current levels.

\begin{figure}[!t]
\centering
\includegraphics[width=\linewidth]{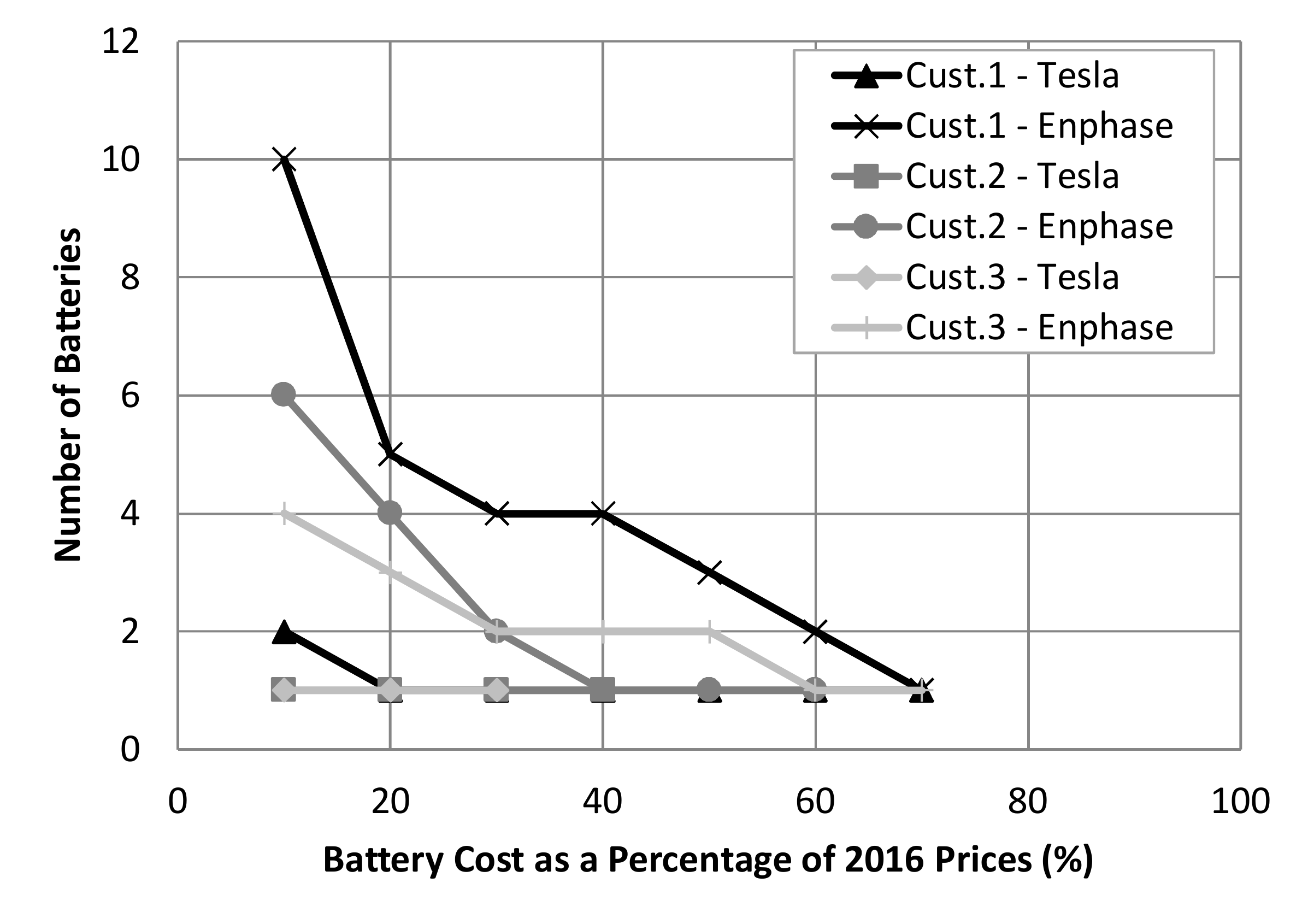}
\caption{Number of batteries in PV-battery system for varying installed battery costs (Retailer B, Mode 2 operation)}
\label{fig:battnum}
\end{figure}

Fig.~\ref{fig:battnum} shows the number of batteries that constitute the optimal system as battery prices are decreased. For the more modular Enphase battery system, an increase in battery quantity is observed for each customer as prices decrease. Customer~1, having a relatively high energy demand, would benefit the most from a larger number of Enphase batteries at each price point compared to the other two customers. In contrast Customer 3 would not benefit from additional batteries until the 30\% price point is reached, after which the customer could can take immediate advantage of additional units as prices continue to decrease. However for the larger Tesla system, a single battery was found to be sufficient for all customers under all battery price scenarios with one exception being an additional battery for Customer 1 at the 10\% price point.

\begin{table}[!t]
\renewcommand{\arraystretch}{1.3}
\caption{Economic Performance Under Different Battery Operating Modes (Tesla Batteries, Cost = 10\% of 2016 Prices)}
\label{tab:batmode}
\centering
\begin{tabular}{|c|c|c|c|c|c|c|}
\hline
\multirow{2}{*}{Customer} & \multirow{2}{*}{Op. Mode} & Battery Size & PV Size & \multirow{2}{*}{NPV}\\
& & (kWh) & (kW$_\text{p}$) & \\\hline\hline
\multirow{4}{*}{1} & 1 & 27 & 8.4 & \$17,446\\\hhline{~*{6}{-}}
  & \cellcolor[gray]{0.8}2 & \cellcolor[gray]{0.8}27 & \cellcolor[gray]{0.8}8.4 & \cellcolor[gray]{0.8}\$19,637\\\hhline{~*{6}{-}}
 & 3 & 27 & 7.56 & \$16,721\\\hhline{~*{6}{-}}
 & 4 & 27 & 8.4 & \$18,102\\\hline\hline

\multirow{4}{*}{2} & 1 & 13.5 & 4.2 & \$3.800\\\hhline{~*{6}{-}}
  & \cellcolor[gray]{0.8}2 & \cellcolor[gray]{0.8}13.5 & \cellcolor[gray]{0.8}4.2 & \cellcolor[gray]{0.8}\$4,931\\\hhline{~*{6}{-}}
 & 3 & 13.5 & 4.2 & \$3,269\\\hhline{~*{6}{-}}
 & 4 & 13.5 & 4.2 & \$4,232\\\hline\hline

\multirow{4}{*}{3} & 1 & 13.5 & 4.2 & \$7,313\\\hhline{~*{6}{-}}
  & \cellcolor[gray]{0.8}2 & \cellcolor[gray]{0.8}13.5 & \cellcolor[gray]{0.8}4.2 & \cellcolor[gray]{0.8}\$8,678\\\hhline{~*{6}{-}}
 & 3 & 13.5 & 4.2 & \$6,785\\\hhline{~*{6}{-}}
 & 4 & 13.5 & 4.2 & \$7,652\\\hline
\end{tabular}
\end{table}

Table~\ref{tab:batmode} summarizes the effect of battery operation mode on the NPV for each customer. The table considers a significantly deflated battery pricing scenario of 10\% of 2016 installation costs whereby energy arbitrage would yield its greatest benefit among the price points considered in this research. In all instances, Mode 2 was found to produce the highest NPV, i.e. the battery operating to maximize self-consumption of PV generated energy in shoulder and peak periods with no energy arbitrage. Consequently, even with significantly deflated battery costs, under current TOU electricity tariffs and with electricity prices continuing to increase at current rates, a battery system engaging in energy arbitrage was not found to provide any additional economic benefit than a battery system purely used for PV generation load shifting.

\section{Conclusion}

Significant price reductions in PV and battery systems have sparked considerable interest in  hybrid PV-battery solutions at the residential level. However optimal system selection is critical to ensure the economic viability of such systems for a particular customer's energy requirements. 

An optimization tool was developed in this paper and applied to three real-world electricity customers. Based on current PV and battery system prices, no battery system was found to be economically viable for the residences assessed, however optimized PV-only systems were found to yield a net benefit for all customers.

A sensitivity analysis was conducted on battery pricing to determine the price point at which a hybrid PV-battery system would yield a net benefit improvement. The results showed that significant price reductions to 60-70\% of current prices are required before the tested customers could take advantage of an energy storage system. It was also concluded that customers can generally take advantage of a modular system of smaller batteries earlier than a bulk energy storage system. Additionally, the results indicated that the current size of the Tesla Powerwall 2 battery is large enough for most energy storage needs even with battery prices at significantly deflated levels.

Finally, various battery operating modes were examined to determine the most economically beneficial operation. No instances were found whereby energy arbitrage yielded a greater benefit than purely maximizing PV self-consumption. This observation continued to hold at all battery pricing levels.





%

\bibliographystyle{IEEEtran}
\bibliography{IEEEabrv,iciea2017}

\end{document}